\documentclass[invmat]{svjour}
\usepackage{latexsym}
\usepackage{graphics}
\usepackage{amsmath}
\usepackage{amssymb}
\usepackage[normalem]{ulem}
\usepackage{amscd}
\usepackage{times}
\begin{document}
\newcommand{\tit}[5]{{\sl #5.}\ #1\ {\bf #2}, #3 (#4)}
\newcommand{\CCAL}{\mathcal{C}}
\newcommand{\ECAL}{\mathcal{E}}
\newcommand{\SCAL}{\mathcal{S}}
\newcommand{\MCAL}{\mathcal{M}}
\newcommand{\OCAL}{\mathcal{O}}
\newcommand{\VEC}[1]{\mathbf{#1}}
\newcommand{\TO}{,\ldots,}
\newcommand{\avec}{\VEC{a}}
\newcommand{\bvec}{\VEC{b}}
\newcommand{\cvec}{\VEC{c}}
\newcommand{\rvec}{\VEC{r}}
\newcommand{\Rvec}{\VEC{R}}
\newcommand{\xvec}{\VEC{x}}
\newcommand{\alphavec}{\boldsymbol{\alpha}}
\newcommand{\betavec}{\boldsymbol{\beta}}
\newcommand{\gammavec}{\boldsymbol{\gamma}}
\newcommand{\deltavec}{\boldsymbol{\delta}}
\newcommand{\alphabar}{\overline{\alpha}}
\newcommand{\scal}[2]{\mean{#1,  #2}}
\newcommand{\mean}[1]{\langle #1 \rangle}
\newcommand{\eq}[1]{eq.~(\ref{#1})}
\newcommand{\EQ}[1]{~(\ref{#1})}
\newcommand{\ddist}[1]{|#1|_{\gamma}}
\newcommand{\sect}[1]{Section~\ref{#1}}
\newcommand{\fig}[1]{Fig.~\ref{#1}}
\newcommand{\lem}[1]{Lemma~\ref{#1}}
\newcommand{\cor}[1]{Corrolary~\ref{#1}}
\newcommand{\prop}[1]{Proposition~\ref{#1}}
\newcommand{\conj}[1]{Conjecture~\ref{#1}}
\newcommand{\M}{\ensuremath{\MCAL}}
\newcommand{\half}{\tfrac{1}{2}}
\newcommand{\pred}{ooo}
\newcommand{\mq}[2]{\uwave{#1}{\marginpar{\footnotesize #2}}}
\newcommand{\ww}[1]{\uwave{#1}}
\newcommand{\mm}[1]{\uwave{#1}}
\title{Jamming and geometric representations of graphs}
\author{Werner Krauth\inst{1} \and Martin Loebl\inst{2}}
\institute{CNRS-Laboratoire de Physique Statistique, Ecole Normale Sup\'{e}rieure, Paris
(email: werner.krauth@ens.fr) \and Department of Applied Mathematics and
Institute for Theoretical Computer Science, Charles University, Prague (email:
loebl@kam.mff.cuni.cz)}
%\date{Received: \today / Revised version: date}
\date{Version: 09-June-04}
\maketitle
\begin{abstract}
We expose a relationship between jamming and a
generalization of Tutte's barycentric embedding. This provides a basis
for the systematic treatment of jamming and maximal packing problems on
two-dimensional surfaces.
\end{abstract}
\section{Introduction}
\label{s:intro}
\begin{figure}[htbp]
\begin{center}
\includegraphics{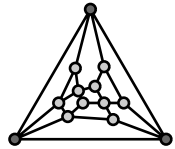}
\end{center}
\caption{Barycentric embedding of a graph with $N=13$ vertices (three outer
vertices).}
\label{f:bary_center_13}
\end{figure}

In a seminal paper \cite{Tutte}, W. T. Tutte addressed the problem of
how to embed a three-connected planar graph in the plane. He proposed
to fix the positions of the vertices on one face and to let the other
(inner) vertices be barycenters of their neighbors
(see \fig{f:bary_center_13}).
The barycentric embedding is unique. It minimizes the energy
\begin{equation*}
E = \sum_{i<j} |\rvec_i - \rvec_j|^2
\label{e:energy_bary}
\end{equation*}
over the positions $\rvec_i$ of the inner vertices.

The purpose of this paper is to study \emph{jamming}, which is
of importance for the physics of granular materials and of glasses
\cite{Torquato,Krauth2000,Krauth2002}, and has many applications in mathematics and computer
science \cite{Conway}. We expose a relationship between jamming  and a
generalization of the barycentric embedding, and provide a basis for the
systematic treatment of jamming on two-dimensional surfaces.

\begin{figure}[htbp]
\begin{center}
\includegraphics{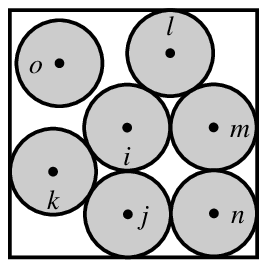} \quad\quad
\includegraphics{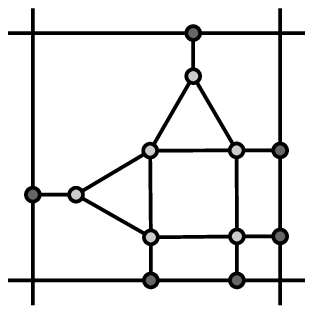}
\end{center}
\caption{Left: Configuration of seven disks in a square.
Disks $i,j,k,l,m,n$ are jammed.
Right: Contact graph of the jammed sub-ensemble. Edges among outer vertices
are omitted.}
\label{f:jammed_config}
\end{figure}

Ensembles of  non-overlapping disks of equal radius may contain
sub-ensembles of disks which do not allow any small moves, regardless of
the positions of the other disks:  In \fig{f:jammed_config} (showing disks
in a square), disks $i,j,k,l,m,n$ are jammed, while $o$ is free to move.

The position $\rvec_i$ of the center of disk $i$ must be more than a disk
diameter away from other disk centers, and more than a disk radius
from the boundary. In a jammed configuration, $\rvec_i$ locally maximizes
the minimum distances to all other disks, and twice the distances to
the boundaries. For a disk $i$ not in contact with the boundary, we have
\begin{equation*}
\min_{j \ne i}| \rvec_i - \rvec_j| = \max_{\rvec} \min_{j \ne i}| \rvec - \rvec_i|.
\label{e:max_min}
\end{equation*}
Equivalently, $\rvec_i$ is a local minimum of the repulsive energy
\begin{equation}
 E(\rvec) = \sum_{j\ne i}| \rvec -\rvec_j|^{q}
\label{e:energy}
\end{equation}
in the limit $q \rightarrow \infty$. Empirical approaches with a repulsive
energy, as in \eq{e:energy}, with very large $q$ have been important to
actually find jammed configurations \cite{Kottwitz}. As the minimum
is local, one cannot prove that with this method all jammed configurations
are generated.

The relationship between jamming and the geometric representation of
graphs was first pointed out, a long time ago, by Sch\"{u}tte and van
der Waerden \cite{Schuette}. Each jammed disk corresponds to an inner
vertex  of a graph, and each touching point with the boundary to an
outer vertex. The edges of the graph refer to contact of disks among
themselves and with the boundary, as shown in \fig{f:jammed_config}.

In this paper, we show that the position $\rvec_i$ is not only the
local maximum of the minimum distance to all other disks, but also the
\emph{global} minimum of the maximum  distance to all the neighbors of $i$
\begin{equation*}
 \max_{a=j,k,l,m} | \rvec_i - \rvec_a| = \min_{\rvec}
 \max_{a=j,k,l,m} | \rvec - \rvec_a|.
\end{equation*}

This leads us to study  \M--representations of graphs (as in
\fig{f:jammed_config}) with inner and possibly outer vertices, where
each one minimizes the maximum (rescaled) neighbor distance rather than
the mean squared distance, as in Tutte's barycentric embedding. Outer
vertices are either fixed or restricted to line segments.

\begin{figure}[htbp]
\begin{center}
\includegraphics{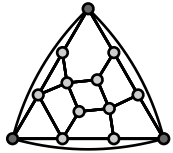}
\end{center}
\caption{Stable \M--representation of the graph of \fig{f:bary_center_13}, with
identical positions of the outer vertices. Three faces of this representation
are flat.}
\label{f:C_infinity_13}
\end{figure}

We define \emph{stable} representations as \M--representations which are local
minima with respect to an ordering relation. This relation replaces the
notion of an energy which cannot be defined in this setting.

We establish that the problem of finding a stable representation in the
plane, torus, or on the hemisphere has an essentially unique solution for
any graph.  We show that \M--representations of three-connected planar and
toroidal graphs are convex pseudo-embeddings, and that the set of regular
three-connected stable representations contains all jammed configurations.
This puts jamming in direct analogy with the barycentric embeddings.
On the sphere, stable representations are not unique, but we conjecture
that their structure is restricted.

One application of jamming is the generation of packings of $N$
non-overlapping disks with maximum radius.  Such a maximal packing
contains  a non-trivial jammed sub-ensemble, since otherwise we could
increase the radius of each disk.  The remaining disks of a maximal
packing are not jammed (as disk $o$ in \fig{f:jammed_config}), and
confined to holes in the jammed sub-ensemble.  In these holes, we can
again search for jammed configurations with suitably rescaled radii.
This gives a recursive procedure to compute maximal disk packings,
which relies on the enumeration of (three--connected) planar or toroidal
graphs and a computation of their jammed representations which form a
subset of the stable representations. Practically, we generate the stable
representations with a variant of the minover algorithm \cite{Krauth1987}
which appears to always converge to a stable solution, on the plane,
torus, and on the sphere.

The most notorious instance of maximal packing is the $N=13$ spheres
problem for disks on the sphere. It has been known since the work of
Sch\"{u}tte and van der Waerden \cite{Schuette} that 13 unit spheres
cannot be packed onto the surface of the unit central sphere (a popular
description has appeared recently in the French edition of Scientific American
\cite{Berger}). However, the minimum radius of the central sphere
admitting such a packing is still unknown, as it is for all larger $N$,
with the exception of $N=24$ \cite{Robinson}.  The problem  of packing
spheres on a central sphere is clearly equivalent to the problem of
packing disks on a sphere.

Our strategy for solving the maximal packing problem will be complete
once the following conjectures are validated:
\begin{conjecture}
\label{c:algorithm}
There exists a finite algorithm to find a stable  representation of a
given graph.
\end{conjecture}
\begin{conjecture}
\label{c:sphere_unique}
Each graph on the sphere with a fixed set of edges crossing a given
equator has at most one non-trivial stable representation up to symmetry
transformations on the sphere.  Furthermore, jammed configurations
are stable.
\end{conjecture}
At present, we are able to prove \conj{c:sphere_unique} for a fixed
\emph{representation} of edges across an equator, rather than their
set. The conjecture is backed by extensive computational experiments.
For planar region and torus, only \conj{c:algorithm} is needed.

\section{\M--representations}
\label{s:M_representations}

In this section we discuss representations of graphs in a planar region,
on the torus, and the sphere. By \emph{torus} we mean a rectangular
planar region where the parallel sides are formally identified.  In a
representation, each vertex is a point, and each edge the shortest
connection between vertices. On the torus, the rectangle can always be
chosen so that vertices  do not lie on its sides.  We require that in a
representation the shortest connection between vertices connected by an
edge is uniquely determined. A representation is  an \emph{embedding}
if it corresponds to a proper drawing.  We recall that a graph is
$k$-connected if it has more than $k$ vertices and remains connected
after deletion of any subset of $k-1$ vertices.  Furthermore, we use
two basic facts of graph-theory: the faces of a two-connected planar
embedding are bounded by cycles, and embeddings of a three-connected
planar graph have a unique list of faces and incidence relations.

\begin{figure}[htbp]
\begin{center}
\includegraphics{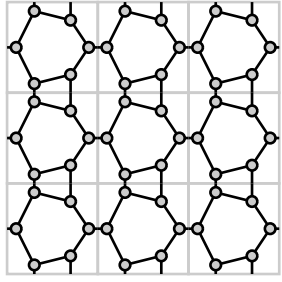}
\end{center}
\caption{The periodic representation of a toroidal graph with six vertices.}
\label{f:periodic_graph}
\end{figure}

Each toroidal representation of a graph gives rise to a unique periodic
representation by tiling the plane with the rectangles, as shown in
\fig{f:periodic_graph}. A \emph{proper} toroidal representation has
edges crossing each side of the rectangle and no outer vertices.

\begin{definition}[Inner and outer vertices] Each graph contains a
possibly empty subset $\OCAL$ of \emph{outer} vertices.  Each outer vertex
is fixed or constrained to lie on a line segment such that any choice of
outer vertices forms a subdivision of a convex $n$--gon, $n\leq |\OCAL|$.
\end{definition}

As indicated in the introduction, we define a rescaled distance, in
order to treat outer and inner vertices on the same footing.

\begin{definition}[Rescaled distance]
The distance between vertices $i$ and $j$ is
\begin{equation*}
\ddist{\rvec_i - \rvec_j} = \gamma_{ij} |\rvec_i - \rvec_j|,
\end{equation*}
\end{definition}
where $\gamma_{ij}$ are positive constants, and $|\quad |$ denotes the
Euclidean distance. In a given representation, we denote by
$l(e)$ the (rescaled) length of an edge $e$, i.e. the distance
between its end-vertices.

\begin{figure}[htbp]
\begin{center}
\includegraphics{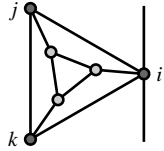}
\end{center}
\caption{Representation of planar graph in
the plane. The outer vertex $i$ is constrained to lie on a line segment,
whereas $j$ and $k$ are fixed.}
\label{f:planar_graph_inner_outer}
\end{figure}

\subsection{\M--center of vectors}
\label{s:M_center}
\begin{definition}[Radius]
The \M--center of a finite number of vectors $\rvec_i,\ i \in I$ is the
vector $\rvec_{\infty}$ minimizing the \emph{radius} of $\rvec$:
$\rho(\rvec) = \max_{i} \ddist{\rvec - \rvec_i}$
\begin{equation*}
\rho(\rvec_{\infty}) = \min_{\rvec}\ \rho(\rvec).
\label{d:centers}
\end{equation*}
\end{definition}
Note that the \M--center is uniquely determined: If there were two
\M--centers with the same radius $\rho$, then the intersection of the
corresponding circles of radius $\rho$ would contain all the neighbors,
but this  intersection is contained in a circle of smaller radius.

\begin{lemma}[No local minimum besides global one]
\label{l:no_local}
If $\rvec$ is not the \M--center of vectors $\rvec_i,\ i \in I$, then
for each $\delta > 0$ there is a vector $\rvec'$ with $|\rvec' - \rvec|
< \delta$ such that $\rho(\rvec) >  \rho(\rvec')$.
\label{l:cinf_center}
\end{lemma}
To determine the \M--center of vectors $\rvec_i,\ i\in I$, we may
consider all pairs and triples of vectors, construct the (unique) circle
passing through them and check that no vertex is outside the circle.
The center of the smallest circle such obtained is the \M--center.

\begin{definition}[\M--representation]
\label{d:m_representation}
An \M--representation of a graph is a representation where each
inner vertex is the \M--center of its neighbors.
\end{definition}

\begin{definition}[Pseudo--embedding]
\label{d:pseudo_embedding}
A representation of a graph is called {\em pseudo--embedding} if it
is an embedding except that some faces may collapse into a line.
Such faces will be called \emph{flat}.
Moreover, a \emph{convex pseudo--embedding} has convex faces and each flat
face is a topological subdivision of $C_2$, a cycle of length two.
\end{definition}
An example of a convex pseudo--embedding is shown in \fig{f:C_infinity_13}.

\begin{proposition}
\label{p:three_conn_embed}
Let $\ECAL$ be a representation of a three-connected planar graph on a
plane or hemisphere, such that each inner vertex belongs to the
convex hull of its neighbors, with non-empty set of outer vertices. Then
$\ECAL$ is a convex pseudo--embedding.
\end{proposition}

\begin{proof}
We proceed analogously to the paragraphs 6-9 of \cite{Tutte}. Since
$G$ is three-connected, its set of faces is uniquely determined, and each
face is bounded by a cycle. $\OCAL$ denotes the set of outer vertices.

Let $l$ be a line in the plane or a non-trivial intersection of 
a plane with the hemisphere and define $g(v), v\in V$, 
as the perpendicular distance of $v$ to $l$, 
counted positive on one side and negative on the other side of $l$.

The outer vertices with the greatest value of $g$ are called {\em positive
poles} and those with the least value of $g$ are {\em negative poles}.
Note that the sets of positive and negative poles are disjoint since $\OCAL$
forms a subdivision of a convex $n$--gon.

A simple path $P=v_1,\dots, v_k$ of $G$ is {\em right (left) rising} if
for each $i$, $g(v_i)< g(v_{i+1})$ or $g(v_i)= g(v_{i+1})$ and $v_{i+1}$
is on the right (left) hand-side of $v_i$.
Right (left) falling paths are defined analogously.

\begin{lemma}
\label{l:vertex_different}
Each vertex $v$ of $G$ different from a pole has two neighbors $v'$
and $v''$ so that $g(v')< g(v)< g(v'')$ or $g(v')=g(v)=g(v'')$ and $v$
belongs to the line between $v'$ and $v''$.
\end{lemma}

\begin{proof}
This follows for outer vertices since they form a subdivision of a convex
$n$--gon, and for inner vertices because of the convexity assumption.
\end{proof}

\begin{lemma}
\label{l:right_left_rising}
Let $v$ be a vertex of $G$. There is a right rising and a left
rising path from $v$ to a positive pole, and also both right and left
falling paths from $v$ to a negative pole.
\end{lemma}

\begin{proof}
By \lem{l:vertex_different} $v$ has a neighbor $v'$ with $g(v')> g(v)$
or $g(v')= g(v)$ and $v'$ is on the right hand-side of $v$. Since $G$
is three-connected, $v'$ has a neighbor different from $v$. Using
\lem{l:vertex_different}, we can monotonically continue from $v'$.
This constructs a right rising path, and the remaining
paths may be obtained analogously.
\end{proof}

\begin{lemma}
\label{l:convex_hull}
If $v\notin \OCAL$ then $v$ belongs to the convex hull of $\OCAL$.
\end{lemma}

\begin{proof}
If such $v$ does not belong to the convex hull of $\OCAL$, then let $l$ be a
line in the plane (cycle on the hemisphere) which defines a separating plane,
 and we get a contradiction with \lem{l:right_left_rising}.
\end{proof}

\begin{lemma}
\label{l:no_cut}
Let $F$ be a face of $G$ and $v_1,v_1',v_2,v_2'$ vertices of $F$ appearing
along $F$ in this order. Then $G$ does not have two disjoint 
$v_1,v_2$ and $v_1',v_2'$ paths.
\end{lemma}
\begin{proof}
This is a simple property of a face of a planar graph.
\end{proof}

\begin{lemma}
\label{l:flat_subdivision}
If a face $F$ is flat then it is a topological subdivision of $C_2$.
Furthermore, let $e$ be an edge of a face $F$ and let $l$ be a line 
in the plane (a cycle on the hemisphere) containing
$e$. Then $F$ is embedded on one side of $l$.
\end{lemma}
\begin{proof}
This simply follows from \lem{l:right_left_rising} and \lem{l:no_cut} (see \fig{f:rising_paths}).
\end{proof}

\begin{figure}[htbp]
\begin{center}
\includegraphics{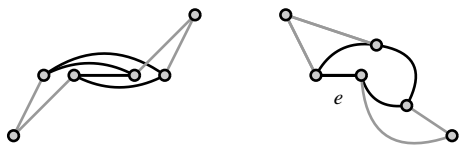}
\end{center}
\caption{Left: a flat face must be a subdivision of $C_2$. Right: each
face must lie on one side of incident edge $e$.}
\label{f:rising_paths}
\end{figure}

It follows from \lem{l:flat_subdivision}  that each face is a subdivision of a convex
$n$--gon or flat and subdivision of $C_2$.

Each edge belongs to exactly two different faces. An edge is
\emph{redundant} if it belongs to two flat faces. More generally, 
in a two-connected representation with prescribed faces such that 
each flat face is a subdivision of $C_2$, a path which is a
subdivision of an edge is {\emph redundant} if it belongs to two flat faces.
A graph is a \emph{simplification} of $G$ if some redundant edges and,
thereafter, maximal redundant paths have been deleted.

\begin{lemma}
\label{l:simplification_c_2}
A flat face of a simplification of $G$ is a subdivision of $C_2$.
\end{lemma}

\begin{proof}
If we delete $e$ and unify the two faces containing $e$, we get a
planar graph.  If the statement does not hold then we can again use
\lem{l:right_left_rising} and \lem{l:no_cut} to obtain a contradiction. The same
applies for a maximal redundant path.
\end{proof}

Let $G'$ be the smallest simplification of $G$ and let $F$ be a flat face
of $G'$. We know that it is a subdivision of $C_2$, and each edge of $F$
belongs to one of the two sides  of $C_2$.

\begin{lemma}
\label{l:claim10}
Let $e$ be an edge of $G'$  and let $l$ be the line in the plane
(cycle on the hemisphere) containing $e$.
\begin{enumerate}
\item[1.]
If $e$ belongs to a flat face $F$, then the faces incident with
edges of different sides of $F$ are on opposite sides of $l$.
\item[2.]
If edge $e$ does not belong to a flat face, then the two faces
incident with $e$ lie on opposite sides of $l$.
\end{enumerate}
\end{lemma}
\begin{proof}
For the second property: as in the proof of \lem{l:simplification_c_2}, if we delete $e$
and 'unify' the two faces containing $e$, we get a planar graph. If the
two faces lie on the same side of $e$, we can use \lem{l:no_cut}. The first
property is analogous.
\end{proof}

Let $|G|$ denote the subset of the surface consisting of the embeddings of
the vertices and edges of $G$, and let $S$ denote the complement of $|G|$.

We define a function $d$ on $S$ as follows: $d(x)=1$ if $x$ is not within
the convex hull of $\OCAL$, otherwise, $d(x)$ equals the number of
interiors of faces to which $x$ belongs. The correctness of this definition
is guaranteed by \lem{l:convex_hull}.

\begin{lemma}
\label{l:d_equals_1}
For each $x\in S$, $d(x)=1$.
\end{lemma}
\begin{proof}
It follows from \lem{l:claim10} that the function $d$ does not change when
passing an edge. However, it cannot change elsewhere and outside of
the convex hull of $\OCAL$ it equals to $1$. Hence it is $1$ everywhere.
\end{proof}

\begin{lemma}
\label{l:edge_interior}
If an edge $e$ intersects the interior of an edge $e'$, then one of
them is not in $G'$ or they belong to opposite sides of a flat face of $G'$.
\end{lemma}
\begin{proof}
This is a corollary of \lem{l:d_equals_1}.
\end{proof}
$\Box$
\end{proof}

\begin{corollary}[\M--rep. is pseudo--embedding]
\label{c:m_representations_pseudo}
An $\M$--representation without outer vertices of three-connected planar
graphs on a sphere or three-connected proper toroidal graphs on a torus
is a convex pseudo--embedding.
\end{corollary}

\begin{proof}
As the number of vertices is finite, we can always find a cut (rectangle
and plane through the center, respectively), which does not contain any
intersection of two edges. The corollary follows by taking as outer
vertices the intersection of edges with the cut. If the cut does not
intersect any edges, the representation is trivial.
\end{proof}

\begin{definition}[Ordering of representations]
\label{d:ordering}
Consider two representations $\ECAL$ and $\ECAL'$ of a graph $G$.
We say that $\ECAL$ is \emph{smaller} than $\ECAL'$ ($\ECAL < \ECAL'$)
if the ordered vector of lengths of the edges containing an inner vertex
of $\ECAL$ is lexicographically smaller than the ordered vector of
lengths of the same edges in $\ECAL'$.
\end{definition}

The above ordering relation cannot generally be mapped into the real
numbers, because the real axis does not admit an uncountable number
of disjoint intervals.
Therefore, there is no `energy' (generalizing
\eq{e:energy_bary}) such that $\ECAL < \ECAL' \Leftrightarrow E(\ECAL)
< E(\ECAL')$.

\begin{definition}[Stable representation]
\label{d:minimum}
Consider a representation $\ECAL$ of a graph $G=(V,E)$ with inner, and possibly
outer vertices $i$ at positions $\rvec_i$. $\ECAL$ is \emph{stable}
if there exists a value $\delta$ such that all embeddings $\ECAL'$ of $G$
with vertices at $\rvec_i'$ with $|\rvec_i -\rvec'_i| < \delta\ \forall i$
satisfy $\ECAL' \ge \ECAL$.
\end{definition}

\begin{proposition}
Stable representations are \M--representations.
\end{proposition}

\begin{proof}
Let the vertex $i$ of $\ECAL$ have the radius $\rho_i$ and let edge
$\{i,j\}$ have length $\rho_i$.  Note that $\rho_j \ge \rho_i$. If
$i$ is not the \M--center of its neighbors, then it follows from
\lem{l:cinf_center} that there is a representation $\ECAL'$ obtained
from $\ECAL$ by a small move of vertex $i$, such that
$\rho'_i< \rho_i$.  All edges $\{k,l\}$ with length bigger than $\rho_i$
 are the same in $\ECAL$ and $\ECAL'$.
No edge $\{k,l\}$ of length $\rho_i$ in $\ECAL$ is longer in
$\ECAL'$ and at least one such edge has shortened.
Finally, edges $\{k,l\}$ shorter than $\rho_i$
in $\ECAL$ may become longer in $\ECAL'$. As a result, we have $\ECAL'<
\ECAL$, which is impossible for a stable representation. $\Box$
\end{proof}
\begin{proposition}[Existence of stable representation]
Each graph has a stable representation.
\end{proposition}

\begin{proof}
Let $\delta > 0 $ be a sufficiently small constant.  Define a sequence
of representations $\ECAL_1, \ECAL_2 , \ldots$ as follows: $\ECAL_1$
is arbitrary. If  $\ECAL_i$ is unstable let  $\ECAL_{i+1}$ be a
lexicographically minimal representation where each vertex has moved by
at most $\delta$ (it exists by compactness). In particular
$\ECAL_{i+1}< \ECAL_i$.
Again by compactness, there is a converging subsequence of representations
$\ECAL_j'$ with limit $\ECAL'$.  $\ECAL'$ must be stable since otherwise
for $\ECAL_i'$ very near to $\ECAL'$, there is a closeby representation
$\overline{\ECAL}< \ECAL'$. Taking into account the minimality rule in
the construction of the sequence of representations, this  contradicts
the assumption that $\ECAL_{i}'$ monotonically decreases in lexicographic
order to $\ECAL'$. $\Box$
\end{proof}

\subsection{Uniqueness of stable representations}
\label{s:uniqueness}

\prop{p:three_conn_embed} implies that each \M--representation is
a convex pseudo--embedding.  Whereas Tutte's barycentric embedding is
unique, the \M--representations are not necessarily unique, as can be
seen by the counter-example sketched in \fig{f:planar_graph_inner_outer}
(the vertices of the inner triangle are in \M--position; they can be
rotated  and rescaled, to remain in \M--position).  However, there is
a unique  \emph{stable} representation.

\begin{lemma}
\label{l:midpoint_plane}
Consider two-dimensional vectors $\rvec_1, \rvec_2,\rvec'_1, \rvec'_2$
with
\begin{align*}
\rvec_1=(x_1,y_1),\ \text{etc}
\end{align*}
and two midpoints $\overline{\rvec}_1$ and $\overline{\rvec}_2$:
\begin{align*}
\overline{x}_1= \half (x_1+x'_1)\\
\overline{y}_1= \half (y_1+y'_1).
\end{align*}
We then have
\begin{equation*}
\ddist{\overline{\rvec}_1-\overline{\rvec}_2} \le
\frac{\ddist{\rvec_1-\rvec_2} + \ddist{\rvec'_1-\rvec'_2} }{2}.
\notag
\end{equation*}
We have $\ddist{\rvec_1-\rvec_2}=
\ddist{\rvec'_1-\rvec'_2}=\ddist{\overline{\rvec}_1-\overline{\rvec}_2}$
only for parallel transport: $\rvec_1=\rvec'_1 + \cvec;  \rvec_2= \rvec'_2 + \cvec$.
\end{lemma}

\begin{proof}
Follows  from  triangle inequality $| \avec + \bvec | \le |
\avec | + | \bvec |$, with $\avec=\rvec_1 - \rvec_2$ and $\bvec=
\rvec'_1-\rvec'_2$ with equality only for parallel transport.
\end{proof}
Note that if $|\rvec_1-\rvec_2| \neq |\rvec'_1-\rvec'_2|$, the midpoint
distance $|\overline{\rvec}_1-\overline{\rvec}_2|$ is smaller than $\max_i
|\rvec_i-\rvec_i'|$.

\begin{proposition}[Unique stable representation in the plane]
\label{p:unique_plane}
Each graph $G$ has
a unique stable representation in the plane (up to parallel transport).
\end{proposition}

\begin{proof}
We assume the contrary.  Let representations $\ECAL^0$ and $\ECAL^1$,
realized by vectors $\rvec_i^0$ and $\rvec_i^1$, be two stable
representation.  We can assume $\ECAL^0 \le \ECAL^1 $.

Consider the representations $\ECAL^{\alpha}$ realized by
\begin{equation*}
\rvec_i^{\alpha} = \rvec_i^0 + \alpha \times
\left[ \rvec_i^1-\rvec_i^0 \right]\quad 0 \le \alpha \le 1 .
\label{e:alpha_relation}
\end{equation*}
The representations $\ECAL^{\alpha}$ exist. We denote by $e^0$ and $e^1$ the
representations of edge $e$ in $\ECAL^0$ and $\ECAL^1$, respectively.
Let $e^1_1 \TO e^1_m$  be the ordered vector of edge lengths.  Let $k$
be the smallest index such that $e^0_k$  is not parallely transported
to $e^1_k$.  We observe the following: if $e$ is an edge of $G$ such
that $l(e^1) = l(e^1_k)$, then $l(e^0) \le l(e^1_k)$ since $\ECAL^0 \le
\ECAL^1 $.  It means by \lem{l:midpoint_plane} that
$\ECAL^{\alpha}< \ECAL^1\ \forall \alpha <1$, which implies that
$\ECAL^1 $ is not stable. $\Box$
\end{proof}

\begin{proposition}[Unique stable representation on torus]
\label{p:unique_torus}
Each graph $G$ has
a unique stable representation on the torus if the sets of
edges crossing each boundary are prescribed (up to parallel transport).
\end{proposition}
\begin{proof} The representations $\ECAL^{\alpha}$ of the
previous proof can analogously be applied to the corresponding periodic
representations, both for edges in the inside of one rectangle and for the
edges going across the boundary.
$\Box$
\end{proof}
This means that the number of stable representations of a toroidal graph
is bounded by the number of possible sets of boundary horizontal and
vertical edges.

Next we will discuss uniqueness of stable embeddings on the hemisphere.
The following \lem{l:midpoint_sphere} is a nontrivial variant
of \lem{l:midpoint_plane}: Obviously, the triangle inequality remains
valid in three dimensions, but the midpoint would not lie on the surface
of the sphere.

\begin{lemma}
\label{l:midpoint_sphere}
Consider three-dimensional vectors $\rvec_1, \rvec_2, \rvec_1', \rvec_2'$
on the unit hemisphere with
\begin{align*}
\rvec_1=(x_1,y_1,z_1= + \sqrt{1 - x_1^2 - y_1^2}),\ \text{etc}.
\end{align*}
and two midpoints $\overline{\rvec}_1$ and $\overline{\rvec}_2$  which are defined
with respect to a projection of vectors on the equator $z=0$
\begin{align}
\overline{x}_i=\half (x_i+x_i');\ \overline{y}_i=\half (y_i+y_i');\
\overline{z}_i = \sqrt{1-\overline{x}_i^2 - \overline{y}_i^2}\quad i=1,2.
\label{e:projected_midpoint}
\end{align}
We then have, for the three-dimensional Euclidean distance-squared:
\begin{equation}
|\overline{\rvec}_1-\overline{\rvec}_2|^2 \le
\frac{|\rvec_1-\rvec_2|^2 + |\rvec_1'-\rvec_2'|^2 }{2}.
\label{e:sphere_inequality}
\end{equation}
In\EQ{e:sphere_inequality}, we can have equality only for generalized
parallel transport with $x_1 - x_1'= c (x_2 - x_2')$ and $y_1- y_1' =
c (y_2 - y_2')$ for special values of $c$.
\end{lemma}
\begin{proof}
We can write\EQ{e:sphere_inequality} as
\begin{multline}
\left( \overline{z}_1 - \overline{z}_2 \right)^2
\le -
\left[ \frac{x_1 + x_1'}{2}- \frac{x_2 + x_2'}{2} \right]^2 +
\frac{ (x_1-x_2)^2}{2} + \frac{(x_1'-x_2')^2}{2} \\
 + \text{same terms in $y, y'$}
 + \half(z_1 - z_2)^2 + \half (z_1' - z_2')^2.
\label{e:reworked_inequality_1}
\end{multline}
Explicit calculation shows that the terms on the first row of expression\EQ{e:reworked_inequality_1} is
\begin{multline*}
-
\left[ \frac{x_1 + x_1'}{2}- \frac{x_2 + x_2'}{2} \right]^2 +
\frac{ (x_1-x_2)^2}{2} + \frac{(x_1'-x_2')^2}{2} \\= \tfrac{1}{4}[x_1 - x_2 - x_1' + x_2']^2,
\end{multline*}
which allows to show that\EQ{e:reworked_inequality_1} and \EQ{e:sphere_inequality}
are equivalent to
\begin{multline}
\left( \overline{z}_1 - \overline{z}_2\right)^2
\le
\tfrac{1}{4} \left( x_1 - x_2 - x_1' + x_2' \right)^2 +
\tfrac{1}{4} \left( y_1 - y_2 - y_1' + y_2' \right)^2 \\
+ \tfrac{1}{2} \left( z_1 - z_2\right)^2
+ \tfrac{1}{2}\left( z_1' - z_2'\right)^2.
\label{e:reworked_inequality}
\end{multline}
Furthermore, we have, from \eq{e:projected_midpoint}
\begin{align*}
\overline{z}_i^2 &= 1 - \left(\frac{x_i + x_i'}{2} \right)^2 -
\left(\frac{y_i + y_i'}{2} \right)^2\\
      &=
\tfrac{1}{4} (x_i-x_i')^2 +
\tfrac{1}{4} (y_i-y_i')^2 +
\tfrac{1}{2} (z_i^2+z_i'^2)\quad i=1,2.
\end{align*}
The inequality\EQ{e:reworked_inequality} now follows from the
triangle inequality $(|\avec| - |\bvec|)^2 \le | \avec -
\bvec|^2$ with the four-dimensional vectors
\begin{align*}
\avec&= \left( \frac{z_1}{\sqrt{2}}, \frac{z_1'}{\sqrt{2}},
\frac{x_1-x_1'}{2}, \frac{y_1-y_1'}{2} \right)\\
\bvec &= \left( \frac{z_2}{\sqrt{2}}, \frac{z_2'}{\sqrt{2}},
\frac{x_2-x_2'}{2}, \frac{y_2-y_2'}{2} \right),
\end{align*}
which evidently satisfy $|\avec|= \overline{z}_1$, $|\bvec|= \overline{z}_2$
and $| \avec - \bvec|^2$ equal to the r.h.s. of\EQ{e:reworked_inequality}.
$\Box$
\end{proof}

\begin{proposition}[Unique stable rep. on hemisphere]
\label{p:unique_hemisphere}
Each graph $G$ on the hemisphere with fixed outer vertices has a unique
stable representation.
\end{proposition}
\begin{proof}
Using the definition of midpoints on the sphere from \lem{l:midpoint_sphere}, 
we can define valid 
representations $\ECAL^{\alpha}$ as in the proof of \prop{p:unique_plane}.
Then we still have $\ECAL^{\alpha}< \ECAL^1\ \forall \alpha <1$. It is easy to see
that, with fixed outer vertices, generalized parallel transport is impossible.
$\Box$
\end{proof}

\subsection{Stable representations on the sphere}
\label{s:stable_sphere}

\begin{figure}[htbp]
\begin{center}
\includegraphics{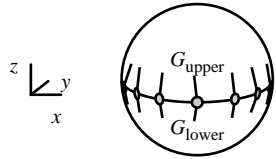}
\end{center}
\caption{An equator embedding of a graph $G$, with a central cycle at $z=0$, and upper and
lower subgraphs $G_{\text{upper}}$ (at $z>0$)  and $G_{\text{lower}}$ (at $z<0$).
The edges on the central cycle are  longer than all other edges.}
\label{f:equator_graph}
\end{figure}

On the sphere, there can be several non-trivial stable representations. To
see this, consider the equator-representation of \fig{f:equator_graph}.
Besides a central cycle, at $z=0$, there are vertices in the upper
subgraph $G_{\text{upper}}$ (with $z>0$) and in the lower subgraph
$G_{\text{lower}}$ (at $z<0$). Furthermore, we suppose that the edges
on the central cycle are longer than those in the rest of the graph. An
equator representation can give rise to two inequivalent representations,
namely by pulling the central cycle \emph{up} to $z>0$, or \emph{down} to $z<0$.

\prop{p:unique_hemisphere} allows us to observe that nevertheless,
some degree of uniqueness can be preserved.

\begin{proposition}[Unique representation with fixed cut]
There is unique stable representation of a graph on a sphere when the
edges crossing a given equator are fixed.
\end{proposition}
As mentioned in the introduction, extensive computing experiments suggest
the stronger statement of \conj{c:sphere_unique}. This would imply that
a given graph has only a finite number of stable representations (up to
symmetry operations).

\section{Jamming}
\label{s:jamming}

In this section we derive basic properties of jammed configurations of
disks on planar region, torus and sphere.

\begin{figure}[htbp]
\begin{center}
\includegraphics{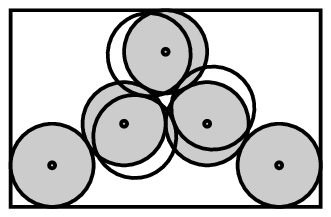}\quad \quad
\includegraphics{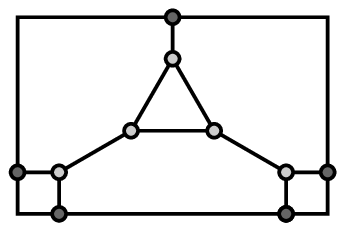}
\end{center}
\caption{Left:  Configuration of five disks, in which no disk can move
by itself, but three disks can move together, as indicated.  Right: The
(unjammed yet stable) embedding corresponding to the configuration. Edges
among outer vertices are omitted.}
\label{f:jammed_rotation}
\end{figure}

\begin{definition}[Jammed embedding] 
\label{prop:jammed} 
An embedding $\ECAL$ of a graph is \emph{jammed} if
\begin{enumerate}
\item $\ECAL$ belongs to the convex hull of the set $\OCAL$ of outer vertices,
and each outer vertex is connected to at least one inner vertex.
\item $\ECAL$ is \emph{regular}, i.e. all edges have
length $\gamma_{ij} \times d$ and the distance between any two vertices $k,l$ not connected by
an edge is strictly bigger than $\gamma_{kl} \times d$.
\item there is  a $\delta$ such that no representation $\ECAL'$ with
$|\rvec_i' - \rvec_i| < \delta$ has some edge longer
and no  edge shorter than in $\ECAL$.
\end{enumerate}
\end{definition}
The three central disks in \fig{f:jammed_rotation} are not jammed,
even though each one cannot move individually.

Let us recall that in a representation, set $\OCAL$ of outer vertices forms
a subdivision of a convex $n$--gon. Let us denote this cycle of outer
vertices by $C_{\OCAL}$ and if $G$ is a jammed embedding then let us denote
by $G_{\OCAL}$ the pseudo--embedding obtained from $G$ by adding the
edges of $C_{\OCAL}$; let us note that $C_{\OCAL}$ bounds the outer
face of $G_{\OCAL}$.
In a jammed embedding $G$ each inner vertex has degree bigger than two 
and at most five on the sphere, and at most six on planar region and torus.  
If $G_{\OCAL}$ is two-connected then each inner face is convex 
(\cite{Schuette}).  
Below we show that each jammed graph $G_{\OCAL}$  is three-connected.

\begin{lemma}
\label{l:two_conn}
If a graph $G$ is connected but not two-connected  and has no vertex of
degree $1$, then there are two vertices $v_1,v_2$ (possibly $v_1=v_2$)
and components $G_i$ of $G \setminus v_i$ ($i=1,2$) such that each
$G_i\cup v_i$ is two-connected and $G_1\cup v_1$ is a subgraph of $G \setminus G_2$.
\end{lemma}

\begin{lemma}
\label{l:conn}
If a graph $G$ is two-connected but not three-connected, and has no vertex of
degree $2$,  then
it has two pairs $V_1, V_2$  of vertices (possibly with $V_1 \cap V_2
\neq \emptyset$) and components $G_i$ of $G \setminus V_i$ ($i=1,2$) such that
$G_i\cup V_i$ is two-connected and $G_1\cup V_1$ is a subgraph of $G \setminus G_2$.
\end{lemma}

\begin{definition}
\label{d:convex}
Let graph $G$ be a two-connected regularly embedded graph, let $F$
be a face of $G$ bounded by a cycle $C$ and let $v \in C$.  We say that
$v$ is convex with respect to  $F$ if there is a plane containing $v$
and perpendicular to the surface (planar region, torus, sphere), so that
a small neighbourhood of $v$ in $F$ lies completely in one half-space
defined by the plane.
\end{definition}

We note that if $F$ is a convex region bounded by a cycle then each
vertex of the cycle is convex with respect to $F$.

\begin{lemma}
\label{l:conv1}
Let $F_1, F_2$ be connected regions bounded by cycles $C_1,C_2$ and such that 
none of them covers the whole planar region (torus, sphere, hemisphere),
but $F_1\cup F_2$ do. Then
there are at least three non-convex vertices of $C_1$ with respect to  $F_1$
or of $C_2$ with respect to $F_2$.
\end{lemma}
\begin{proof}
Cycle $C_2$ must be embedded inside $F_1$ and
$F_2$ must be the region defined by $C_2$ that is not a subset of $F_1$. Then
$C_2$ is a non-self-intersecting cycle on $F_1$. As such, it must have at
least three sharp corners on $F_1$.
\end{proof}

\begin{proposition}[Jammed graph three-connected]
\label{p:jammed_three_connected}
If an embedding $G$ is jammed, then $G_{\OCAL}$ is three-connected.
\end{proposition}

\begin{figure}[htbp]
\begin{center}
\includegraphics{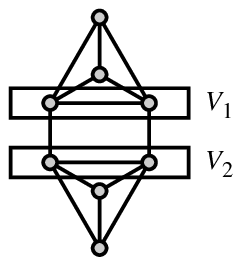}\quad \quad
\includegraphics{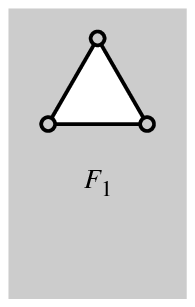}\quad \quad
\includegraphics{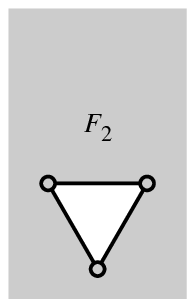}
\end{center}
\caption{ }
\label{f:jammed_graph_three}
\end{figure}

\begin{proof}
Let $G_{\OCAL}$ be a minimum counter-example.
If $G_{\OCAL}$ is not connected then each component is jammed
and three-connected by the minimality assumption. Consider
the embeddings of different components $G_1$ and $G_2$ in the embedding
of $G$. $G_1$ is completely embedded in one of the faces of $G_2$
and vice versa, which is not possible by convexity of faces
and by \lem{l:conv1}.

If $G_{\OCAL}$ has a vertex of degree at most two then it cannot be
jammed. Therefore, we suppose that $G_{\OCAL}$ is without a vertex of degree
two and either connected or two-connected. By \lem{l:two_conn} and
\lem{l:conn} there is a subset of vertices $V_1$
and a component $G_1$ of $G_{\OCAL} \setminus V_1$ such that $G_1\cup
V_1$ is two-connected and has one of the following two properties:

\noindent
1.$V_1$ consists of a single vertex $v_1$ and  there is a vertex $v_2$ with
$V_2 = \{ v_2\}$
and a component $G_2$ of $G_{\OCAL} \setminus V_2$ so that $G_2 \cup V_2$ is a
two-connected subgraph of $G_{\OCAL} \setminus G_1$.

\noindent
2. $V_1$ consists of two vertices and  there is a subset $V_2$ of
\emph{two} vertices and a component $G_2$ of $G_{\OCAL} \setminus V_2$ so that $G_2\cup V_2$
is a two-connected subgraph of $G_{\OCAL} \setminus G_1$.

We consider the embedding of $G_1\cup V_1$ induced by the embedding
of $G_{\OCAL}$ (see \fig{f:jammed_graph_three}). Let $F_1$ be the face in which
$G_{\OCAL} \setminus G_1$ is embedded and let $C_1$ be the bounding cycle of
$F_1$. Clearly $V_1\subset C_1$.
Moeover $F_1$ cannot be the outer face of the embedding
of $G_{\OCAL}$ since $G_{\OCAL} \setminus G_1$ is embedded there. 
Hence all the vertices of $C_1 \setminus V_1$ are convex 
with respect to $F_1$.

Let $F_2$ be the face of the induced embedding of $G_2\cup V_2$ which
contains $C_1$ and let $C_2$ be the cycle that bounds $F_2$. Clearly $V_2
\subset C_2$ and as above, $F_2$ cannot be the outer face of the embedding
of $G_{\OCAL}$. Hence all vertices of $C_2\setminus V_2$ are convex
with respect to $F_2$. Then $F_1, F_2, C_1, C_2$ satisfy the properties
of \lem{l:conv1} but $V_1$ and $V_2$ have only two vertices,
a contradiction.
$\Box$
\end{proof}

\begin{proposition}[Jammed graphs stable]
\label{p:jammed_stable}
Jammed embeddings are stable on the planar region and torus and on
a hemisphere. On the sphere, jammed embeddings are stable if we fix
representations of edges across an equator.
\end{proposition}

\begin{proof}
We show that an unstable regular embedding
$\ECAL$ has a small move which leaves some edges the same and increases
the others, which means that it is not jammed.
For any $\delta>0$, there is a representation $\ECAL'$ within $\delta$  of
$\ECAL$ such that some edges decrease in length, and the others stay the same.
Let $\ECAL''$ be the inverse of the move which took $\ECAL$ into $\ECAL'$
(on the hemisphere: inverse of the projected move).
The inverse move exists, since a jammed embedding in the plane or on the torus
 cannot have an outer vertex on an end point of the corresponding line segment,
and for the hemisphere, if originally an inner vertex positioned
on the equator moved, then $\ECAL$ could not have been jammed.
For the representation $e$ of an edge in  $\ECAL$, let $e'$ and $e''$ be the
respective representations of the same edge in $\ECAL'$ and $\ECAL''$.
We can use \lem{l:midpoint_plane} (on planar
region and torus) and \lem{l:midpoint_sphere} (for the hemisphere) to show
that
$l(e')=l(e) \implies l(e'') \ge l(e)$ and $l(e')< l(e) \implies l(e'') > l(e)$.
$\Box$
\end{proof}
Note that the converse of \prop{p:jammed_stable} is not true and that
stable representations are not necessarily jammed.
An example is shown in \fig{f:jammed_rotation}.

As mentioned in the introduction, extensive computing experiments suggest,
for the sphere, the stronger statement of \conj{c:sphere_unique}.

\section{Algorithms}
\label{s:algorithm}

In \sect{s:M_center}, we discussed a finite algorithm for the
determination of the \M--center of vectors $\rvec_i$.  This algorithm
is of practical use because a circle  in $d$--dimensional space is
already specified by $d+1$ points.  The incremental \emph{minover} algorithm
\cite{Krauth1987} remains useful in high dimension $d$ and is trivial 
to implement. 
For a finite number of vectors $\rvec_i, i \in I$ on the unit sphere, it is defined by
\begin{equation*}
\Rvec_0=0, \quad \Rvec_{k+1} \leftarrow \Rvec_k + \rvec_{i_{\min}}, 
\end{equation*}
where the index $i_{\min}$ is a minimal overlap (scalar product)
vector with
\begin{equation*}
\scal{\Rvec_k}{\rvec_{i_{\min}}} = \min_{i \in I} \scal{\Rvec_k}{\rvec_{i}}.
\end{equation*}
It can be proven \cite{Krauth1987} that
\begin{equation*}
\Rvec_k/|\Rvec_k| \rightarrow \rvec_{\infty}
\label{e:}
\end{equation*}
under the condition that a vector $\Rvec$ exists with
$\scal{\Rvec}{\rvec_{i}} > 0 \ \forall i\in I$.
\begin{figure}[htbp]
\begin{center}
\includegraphics{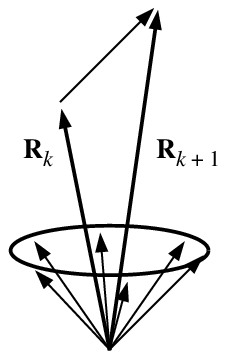} 
\quad \quad
\quad \quad
\includegraphics{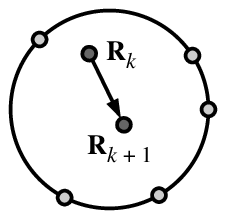}
\end{center}
\caption{Left: Minimum overlap algorithm. 
The rescaled vector $\Rvec_k/|\Rvec_k|$ converges to the \M--center of vectors on the 
sphere.  Right: On the plane, the move is in direction of the vector $\rvec_i$
with maximum distance to $\Rvec_k$. The amplitude $\epsilon_k$ of the move
decreases with $k$, but $\sum_k \epsilon_k$ diverges.}
\label{f:minover_sphere}
\end{figure}

In the  minover algorithm  (on the sphere),  the corrections to $\Rvec_k$ 
decrease with increasing $k$. 
On the plane or the torus, we can do the same by using an update 
\begin{equation*}
\Rvec_{k+1} \leftarrow \Rvec_k + \epsilon_k\   [\rvec_{i_{\max}}- \Rvec_k],
\label{e:max_dist}
\end{equation*}
where $\rvec_{i_{\max}}$ is the vector of maximum distance to $\Rvec_k$ (see
\fig{f:minover_sphere}), and where the sequence $\epsilon_k$ satisfies 
the conditions:
\begin{gather*}
\epsilon_k \rightarrow 0 \ \text{for $k \rightarrow \infty$}\\
\sum_k \epsilon_k \rightarrow \infty \ \text{for $k \rightarrow \infty$}.
\end{gather*}
This algorithm was applied to all vertices sequentially in order to 
compute stable \M--representations.

\begin{figure}[htbp]
\begin{center}
\includegraphics{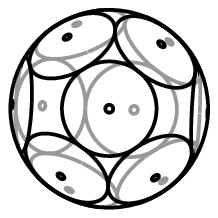}
\end{center}
\caption{Left:  Conjectured optimal configuration of $13$ disks on a unit sphere,
and corresponding representation of vertices, obtained by computational
experiment as stable \M--representation of the graph of \fig{f:bary_center_13} and
\fig{f:C_infinity_13}.  The algorithm of \sect{s:algorithm} was used.}
\label{f:13_disks}
\end{figure}

As a simple test, we have run this algorithm on the three-connected graph
of \fig{f:bary_center_13} and \fig{f:C_infinity_13}, embedded on a sphere
and starting from an equator position. This graph, incidentally,
corresponds to the conjectured optimal packing for the thirteen-sphere problem.
The algorithm converges rapidly to the conjectured optimum solution
\cite{Kottwitz}, which is shown in \fig{f:13_disks}.

As stated in \conj{c:algorithm}, we are convinced that a finite algorithm
for computing a stable \M--representation exists.

\begin{acknowledgement}
We are indebted to F.~Bouchut for the proof of
\lem{l:midpoint_sphere}  and to \mbox{J.~Mlcek} for discussion on
lexicographic order and energy.  This work has been supported by
the European Science Foundation (ESF) program 'Phase Transitions and
Fluctuation Phenomena for Random Dynamics in Spatially Extended Systems'.
We thank ITI Prague and LPS-ENS for hospitality.

\end{acknowledgement}

\end{document}